\documentclass[a4paper,twoside]{article}
\usepackage{amsmath}
\usepackage{amsfonts}
\usepackage{amssymb}
\usepackage{graphicx}
\usepackage{epstopdf}
\usepackage{algorithmic}
\ifpdf
  \DeclareGraphicsExtensions{.eps,.pdf,.png,.jpg}
\else
  \DeclareGraphicsExtensions{.eps}
\fi

\newcommand{\Jac}{\mathrm{Jac}}
\newcommand{\Sing}{\mathrm{Sing}}

\newcommand{\Q}{\mathbb{Q}}
\newcommand{\R}{\mathbb{R}}
\newcommand{\SE}{\mathrm{SE}}
\newcommand{\SO}{\mathrm{SO}}
\newcommand{\T}[1]{#1^{\mathsf{T}}}
\newcommand{\Proj}{\mathbb{P}}

\newtheorem{theorem}{Theorem}
\newtheorem{proposition}[theorem]{Proposition}
\newtheorem{corollary}[theorem]{Corollary}
\newtheorem{lemma}[theorem]{Lemma}

\title{Rationality of the locus of singularities of the general Gough-Stewart platform}
\author{Michel Coste\thanks{Univ Rennes, CNRS, IRMAR - UMR 6625, F-35000 Rennes, France} and Seydou Moussa\thanks{Univ. Dan Dicko Dankoulodo, D\'ep. Math\'ematiques, Maradi, Niger}}

\begin{document}

\maketitle
\abstract{We prove that the set of singular configurations of a general Gough Stewart platform has a rational parametrization. We introduce a reciprocal twist mapping which, for a general orientation of the platform, realizes the cubic surface of singularities as the blowing up of a quadric surface in five points.


}

\section*{Introduction}

Gough-Stewart platform is the most well-known parallel robot with six degrees of freedom. It consists in a platform linked to a fixed base via six limbs whose lengths are controlled by actuated prismatic joints ; each limb is attached to the base by a universal joint and to the platform by a spherical joint (6-U\underline PS in the standard notation \cite{Mer}). 

\begin{figure}[h]\caption{A Gough-Stewart platform}
\begin{center}
\includegraphics[width=0.5\linewidth]{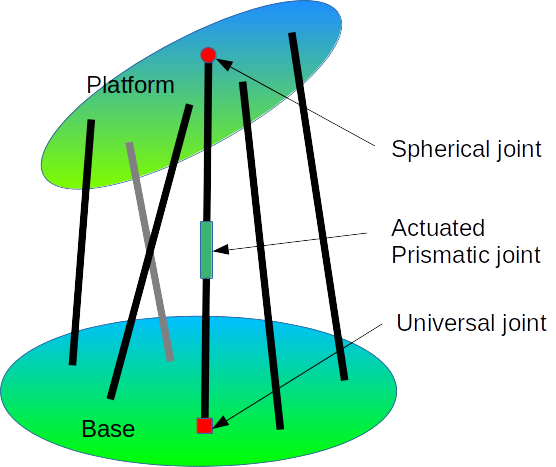}
\end{center}
\end{figure}

The problem of avoidance of singular configurations (in which one loses control on at least one degree of freedom) is crucial for parallel robots, and many papers have been devoted to the study of the set $\Sing$ of singular configurations of the Gough-Stewart platform. A basic fact is the characterization of singular configurations as those for which the systems of Pl\"ucker coordinates of the six limbs of the platform are linearly dependant (see \cite {Mer}). This leads to a description of $\Sing$ by a polynomial equation in terms of position variables and orientation variables. This equation is of degree $3$ in the position variables for a general architecture (this degree may drop for special architectures, see \cite{Ka,Na10}). Hence, for a generic (in the sense of algebraic geometry) orientation, the singular configurations form a cubic surface $\Sigma$ in 3-dimensional space. 

A way to describe the set of singular configurations of a robot is to give a rational parametrization of this set, when this is possible; this has been done for 3-R\underline PR planar robots in \cite{HG}, which stresses the advantage of rationality. In the case of the Gough-Stewart platform, a rational parametrization has been found for special architecture in \cite{BG}, and the rationality of $\Sing$ has been proved for the general case of planar base and platform base in \cite{CoMo}. The rationality in this case follows from the existence of a singular point at infinity of the cubic surface $\Sigma$. This argument cannot be used for the general Gough-Stewart platform : the projective closure $\Sigma^h$ of $\Sigma$ is in general a nonsingular cubic surface. It is known (see \cite{SD}) that the rationality of a nonsingular cubic projective surface is related to special features of the set of the 27 lines on this surface. Experiments using exact computations over rational numbers show that the polynomial of degree 27 whose roots corresponds to the lines splits in factors of degrees 2, 5, 10, 10, where the lines corresponding to the factor of degree 2 are skew. This gives an evidence for rationality, which has to be confirmed by a proof. The key to the proof is given by a consideration from kinematics: to each singular configuration of the platform one can associate a line of reciprocal twists, which expresses the infinitesimal rigid motion which can no longer be controlled in the singular configuration. We prove that, for a general architecture, this reciprocal twist mapping is a regular mapping from the projective cubic surface $\Sigma^h$ to a quadric surface $Q$ (in the 5-dimensional projective space of twists). We prove moreover that the reciprocal twist mapping is the blowing-up of the quadric surface $Q$ in five points, and the five exceptional divisor in the cubic surface $\Sigma^h$ are the lines corresponding to the factor of degree 5 mentioned above. The singularities which form these five exceptional divisors are those for which all six limbs of the platform have a common secant line, and the reciprocal twist is a twist of rotation about this secant line.

We have proved in this way the rationality of the locus of singularities of the general Gough-Stewart platform. But, more than that, the information obtained about the reciprocal twist mapping may be of interest for the kinematics of this robot.

The paper is organized as follows. In Section \ref{Sec2} we describe the singularity locus and explain how the problem of its rational parametrization reduces to the problem of the rationality of the cubic surface $\Sigma$ defined over the field of functions of the rotation group $\SO(3)$; we recall the relation between the rationality of this surface and the structure of the set of lines on this surface, and give the example of a computation about these lines which shows an evidence for rationality. In section \ref{Sec3} we introduce the reciprocal twist mapping and show that it gives a birational equivalence between $\Sigma$ and a quadric surface $Q$, which proves the rationality of $\Sigma$. In section \ref{Sec4} we extend the reciprocal twist mapping to the projective closure $\Sigma^h$ of $\Sigma$ and show that the extended mapping is the blowing-up of $Q$ in five points; we make precise the kinematic significance of the exceptional divisors of this blowing-up. We conclude by collecting some open questions raised by our approach.

\section{The singularity locus of a Gough-Stewart platform, and the question of its rational para\-metrization}\label{Sec2}

 We consider a Gough-Stewart platform with an arbitrary architecture. For $i=1,\ldots,6$, we denote by $A_i$ the vector of coordinates of the center of the universal joint of the $i$-th limb on the base, in the fixed orthonormal frame, and by $b_i$ the vector of coordinates of the center of the spherical joint on the mobile platform in the orthonormal frame attached to this platform. We may assume $A_1=b_1=\T{(0,0,0)}$ and use it whenever it is convenient.

\subsection{Equation of the singularity locus in $\SE(3)$}

We recall the derivation of the equation of the singularity hypersurface, mainly in order to fix notation. The computation is very similar to the one in \cite{DGC}.

The group $\SE(3)$ of rigid motions in 3-space acts on the mobile platform by the transformation $X\mapsto R\,X+P$ where $R$ is the rotation matrix and $P=\T{(x, y, z)}$ is the translation vector. So the coordinates of the joints on the mobile platform are, in the fixed frame, $Rb_i+P$ for $i=1\ldots 6$ ; we set $B_i=Rb_i$ and $C_i=B_i-A_i$.

The Pl\"ucker coordinates of the limbs w.r.t. the fixed frame are 6-dimensional vectors whose first three coordinates are $C_i+P$ and last three coordinates (the moment w.r.t. the origin) $A_i \times (C_i+P)$. It is well known and explained in \cite{Mer} that the Gough-Stewart platform is in a singular configuration if and only if the Pl\"ucker coordinates of the six limbs are linearly dependant. This is expressed by the vanishing of the determinant of the $6\times 6$ matrix whose rows are the Pl\"ucker coordinates of the limbs:
\begin{equation}
\Jac=\T{\begin{pmatrix}C_i+P\\ A_i\times (C_i+P)  \end{pmatrix}_{i=1,\ldots,6}}\;.
\end{equation}
We denote by $\Sing\subset \SE(3)$ the hypersurface of singular configurations, whose equation is $\det(\Jac)=0$. The following result is well known.

\begin{theorem} \label{Th1} For a generic Gough-Stewart platform, the singularity locus $\Sing$ in $\SE(3)$ has an equation which is of degree 3 with respect to $P$.
\end{theorem}

\noindent\textit{Proof:}
We compute $\det(\Jac)$ using the generalized Laplace expansion w.r.t. the first three rows of $\Jac$. We use the notation $[U,V,W]$ to denote the mixed product of the three 3-dimensional vectors $U,V,W$ (i.e. the determinant whose columns are $U,V,W$).  
 \begin{equation}\label{jacobien}
\begin{aligned}
\det(\Jac) = &\sum_{1\le i_1< i_2< i_3\le 6} (-1)^{i_1+ i_2+ i_3}  [C_{i_1}+P, C_{i_2}+P, C_{i_3}+P]\\ &\qquad [A_{j_1}\times (C_{j_1}+P), A_{j_2}\times (C_{j_2}+P), A_{j_3}\times (C_{j_3}+P)]\;,
\end{aligned}
\end{equation}
where $j_1<j_2<j_3$ are the integers between $1$ and $6$ different from $i_1,i_2,i_3$. The mixed products $[C_{i_1}+P, C_{i_2}+P, C_{i_3}+P]$ have degree at most $1$ w.r.t. $P$, while the mixed products $[A_{j_1}\times (C_{j_1}+P), A_{j_2}\times (C_{j_2}+P), A_{j_3}\times (C_{j_3}+P)]$ have degree at most $2$, because the vectors $A_{j_1}\times P, A_{j_2}\times P, A_{j_3}\times P$ are linearly dependant. So the degree of $\det(\Jac)$ w.r.t. $P$ is at most $3$. The computation in a example (see later) shows that it is actually $3$ for a generic Gough-Stewart platform.\quad $\square$

\subsection{Cayley parametrization for rotation matrices}\label{ParamCay}

Let $$U=\begin{pmatrix}  0 & -r & q\\  r & 0 & -p\\ -q & p & 0 \end{pmatrix}$$ be a skew-symmetric matrix.
Since $1$ is not an eigenvalue of $U$, $I-U$ is invertible. The matrix $R=(I+U)(I-U)^{-1}$ is a rotation matrix
\begin{equation}R(p,q,r)=\frac{1}{\Delta}\begin{pmatrix}
1+p^2-q^2-r^2 	& 2(p\,q- r) 	& 2(p\,r+q)		\\
2(p\,q+r)		& 1-p^2+q^2-r^2 	& 2(q\,r-p) 		 \\
2(p\,r-q)		& 2(q\,r+p)		& 1-p^2-q^2+r^2
\end{pmatrix}\;,\end{equation}
where $\Delta=1+p^ 2+q^ 2+r^ 2$. If the vector $\T{(p, q, r)}$ is not the zero vector, it spans the axis of the rotation $R$. The tangent of the half-angle of the rotation is $\sqrt{\Delta-1}$. The Cayley parametrization is a rational parametrization of all rotation matrices, except the half-turns. These half-turns are obtained as limits as $p^ 2+q^ 2+r^ 2$ tends to infinity; alternatively, one can use the homogeneous Euler-Rodrigues parametrization (that is, parametrization with quaternions) with one more variable.

Note that we can recover rationally $p,q,r$ from the rotation matrix $R(p,q,r)$. Indeed
\begin{equation}p=\frac{R_{3,2}-R_{2,3}}{1+\mathrm{tr}(R)},\
q=\frac{R_{1,3}-R_{3,1}}{1+\mathrm{tr}(R)},\
r=\frac{R_{2,1}-R_{1,2}}{1+\mathrm{tr}(R)}\;,\end{equation}
where $\mathrm{tr}(R)=R_{1,1}+R_{2,2}+R_{3,3}$. This shows:
\begin{proposition}\label{Caypar}
The Cayley parametrization  induces an isomorphism between the field $\R(\SO(3))$ of rational functions on $\SO(3)$ and the field of rational functions in three independent variables  $\R(p,q,r)$.
\end{proposition}
The preceding result says, in terms of algebraic geometry, that the variety $\SO(3)$ is a rational variety over $\R$.\par\medskip

\subsection{The cubic surface $\Sigma$}

The equation $\det(\Jac)=0$ defines the algebraic variety of singular configurations $\Sing\subset\SE(3)$. For a fixed $R\in \SO(3)$, the equation $\det(\Jac)=0$ in the variables $x,y,z$ defines a cubic surface $\Sing_R\subset \R^3$, and $\Sing$ may be viewed as the family of these cubic surfaces parametrized by $\SO(3)$. We can also view the equation $\det(\Jac)=0$ as an equation in the three variables $x,y,z$ with coefficients in the field $\R(\SO(3))$ (or the field $\R(p,q,r)$, according to Proposition \ref{Caypar}). As such, this is the equation of a cubic surface $\Sigma$ defined over $\R(\SO(3))$. This surface $\Sigma$ will be the main object of study in the rest of this article. In the language of algebraic geometry, $\Sigma$ is the generic fibre of the family of surfaces $\Sing_R$ parametrized by $\SO(3)$.\par
In order to prove that the algebraic variety $\Sing$ is rational over $\R$, it suffices to prove that $\Sigma$ is rational over $\R(\SO(3))$, since $\SO(3)$ itself is rational over $\R$, as shown in Proposition \ref{Caypar}. This will be the aim of the next section. So we shall work over the field $\R(\SO(3))$. We shall feel free to use the classical formulas of vector algebra relating dot product, cross product, mixed product over $\R(\SO(3))$ as if we were in $\R^3$.

A general criterion for rationality of a cubic surface over a field $k$ is due to Swinnerton-Dyer. It is related to the geometry of lines on the cubic surface. A very classical result of algebraic geometry says that a smooth projective cubic surface over a field $k$ has $27$ lines, which are defined over the algebraic closure of $k$ (see for instance \cite{Har}, Theorem 4.9 p.402).

 \begin{theorem}[\cite{SD}]\label{SwDy} Let $\Gamma$ be a smooth projective cubic surface over a field $k$ (of characteristic $0$). Call $S_n$ a set of $n$ of the $27$ lines on $\Gamma$, mutually skew, which is stable under conjugation over $k$; i.e., the union of these lines is defined over $k$. Then $\Gamma$ is $k$-rational iff it has a $k$-rational point and a $S_2$, or a $S_3$, or a $S_6$.\end{theorem}
 
 We shall first investigate experimentally whether the projective closure $\Sigma^h$ of the cubic surface $\Sigma$ may satisfy this criterion. The existence of a $\R(\SO(3))$-rational point on $\Sigma$ is certain, since $P=\T{(0, 0, 0)}$ is clearly always on $\Sigma$ (when we assume $A_1=b_1=\T{(0, 0, 0)}$). So the problem is only the splitting of the set of 27 lines over $\R(\SO(3))$.
  
 In principle, we could do the computation leaving the geometry of the Gough-Stewart platform free and the orientation $R$ free. This would be a too heavy computation. Instead, we fix the geometry of the platform by choosing rational coordinates and pick a rotation matrix $R$ also with rational coefficients, so that all the computations will be exact computations over $\Q$, concerning the surface $\Sing_R$ defined over $\Q$. If the set of lines on $\Sigma$ splits over $\R(\SO(3)$ then the set of lines on $\Sing_R$ will split accordingly over $\Q$. Experiments show the presence of a $S_2$ defined over $\Q$ on $\Sing_R$, with the same pattern of splitting for different choices of $R$ and different geometries of the platform. Note that for a general cubic surface defined over $\Q$, the set of 27 lines does not split over $\Q$.
 
 \subsection{A case study}\label{Case}
 We fix the geometry of the Gough-Stewart platform as follows:
$$A_1=\begin{pmatrix}0\\0\\0\end{pmatrix} 
A_2=\begin{pmatrix}2\\0\\0\end{pmatrix}
A_3=\begin{pmatrix}0\\2\\-1\end{pmatrix}
A_4=\begin{pmatrix}0\\1\\2\end{pmatrix}
A_5=\begin{pmatrix}1\\0\\1\end{pmatrix}
A_6=\begin{pmatrix}6\\3\\0\end{pmatrix}$$
$$b_1=\begin{pmatrix}0\\0\\0\end{pmatrix} 
b_2=\begin{pmatrix}2\\3\\-1\end{pmatrix}
b_3=\begin{pmatrix}0\\1\\4\end{pmatrix}
b_4=\begin{pmatrix}1\\3\\1\end{pmatrix}
b_5=\begin{pmatrix}1\\3\\-1\end{pmatrix}
b_6=\begin{pmatrix}2\\4\\-3\end{pmatrix}$$
We also choose Cayley parameters $p=0,\;q=0,\;r=0$ for the rotation matrix $R(p,q,r)$, so that $R$ is the identity matrix and $B_i=b_i$ for $i=1,\ldots,6$ 

The set of singular configurations for the chosen rotation $R$ is the cubic surface $\Sing_R$ with equation
\begin{equation}
\begin{aligned}
&80x^3-107yx^2-47zx^2-376x^2-9y^2x-301yx+95zxy-1392x \\
&\quad{} -96z^2x+643zx-68y^3+98zy^2-78y^2+426y+708zy\\
&\quad{}+78z^2y+234z^2-24z^3+1410z
\end{aligned}
\end{equation}
In order to check that the projective closure $(\Sing_R)^h$ of this cubic surface is nonsingular, we homogenize the cubic equation with homogeneization variable $w$ and verify that the ideal generated by the partial derivatives of the homogenized equation w.r.t. $w,x,y,z$ contains the fifth power of the  ideal generated by $w,x,y,z$.

There are 27 lines on this cubic surface; it happens that they are all real for the present example, but this is not the case for other choices. These lines can be computed in the following way: we write a parametrization of a line with parameter $t$ and indeterminate coefficients $a,b,c,d$
$$ x=t,\ y=a+b\,t,\ z=c+d\,t\;.$$
Substituting this parametrization in the equation of the cubic surface yields a degree 3 polynomial in $t$ whose coefficients depend on $a,b,c,d$. The ideal generated by these four coefficients is the ideal of lines on the surface (assuming that no line is contained in the plane at infinity nor parallel to the $(y,z)$-plane). We check that the quotient of $\Q[a,b,c,d]$ by this ideal is indeed a finite extension of degree 27. The Groebner basis of the ideal with respect to the lexicographic order on $a,b,c,d$ contains a polynomial of degree 27 in $d$, and $a,b,c$ are rational polynomials in $d$ in the quotient. This polynomial is the product of four irreducible factors over $\Q$, of degrees respectively 2, 5, 10 and 10. The factors of degree 2 and 5 are
\begin{equation}\label{27droites}
\begin{aligned}
F_2&=2796d^2+4137d-56\\
F_5&=14853594d^5+160133255d^4-6870509d^3\\
&\qquad {} -1145865348d^2+1491086416d-515006656
\end{aligned}
\end{equation}
We explore the coplanarity relations between the lines corresponding to the different factors. In order to do this, we make computations over the quadratic extension of $\Q$ given by $F_2$ and we use the fact that the lines with coeffcients $a,b,c,d$ and $a',b',c',d'$ are coplanar if and only if 
$$(a-a')(d-d') -(b-b')(c-c')=0\;.$$
It can be checked that 
\begin{itemize}
\item the two lines corresponding to the factor $F_2$ are skew, 
\item the five lines corresponding to the factor $F_5$ are those among the 27 which meet both lines of $F_2$, 
\item the ten lines corresponding to one of the factor of degree $10$ are those which meet exactly one of the lines of $F_2$ (five each),
\item the ten lines corresponding to the remaining factor are those which meet none of the lines of $F_2$.
\end{itemize}
\begin{figure}[h]\caption{The cubic surface $\Sing_R$, the two lines corresponding to $F_2$ in black, the five lines corresponding to $F_5$ in red}
\begin{center}
\includegraphics[width=0.7\linewidth]{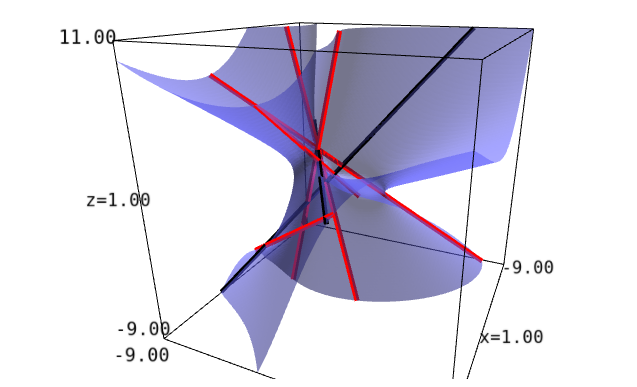}
\end{center}
\end{figure}

The same factorization pattern appears for generic choices of architecture and orientation (with rational parameters).
The two lines corresponding to factor $F_2$ form a $S_2$ over $\Q$. This is a strong indication that $\Sigma$ itself has a $S_2$ over $\R(\SO(3))$, and hence is rational over this field. We shall prove this fact in the following section.

It is well-known how to obtain a rational parametrization from a $S_2$. Pick a point on the cubic surface ; there is a unique line through this point meeting the two lines of the $S_2$, intersection of the planes containing the point and a line of the $S_2$. The set of lines meeting the two lines of the $S_2$ can be seen as a quadric surface, and we obtain in this way a regular mapping from the cubic surface to a quadric surface; this regular mapping is a birational isomorphism, actually the blowing-up of the quadric surface in the five points corresponding to the five lines contained in the cubic surface and meeting the two lines of the $S_2$. For more detail on this, see for instance \cite{Mum}. We shall in the next section show how kinematics provides a regular birational morphism from $\Sigma$ to a quadric surface, and we shall prove in Section \ref{Sec3} That the projective closure $\Sigma^h$ is indeed the blowing up of this quadric surface in five points.

\section{The reciprocal twist at a singular configuration}\label{Sec3}

\subsection{Solutions of homogeneous systems of corank 1}

Recall that a twist $\begin{pmatrix} \Omega\\ V\end{pmatrix}$ is said to be reciprocal to a screw $\begin{pmatrix} F\\ M\end{pmatrix}$ when their reciprocal product $F\cdot V+ \Omega\cdot M$ is zero. If a Gough Stewart platform is in a singular configuration such that the rank of the family of systems of Pl\"ucker coordinates $\begin{pmatrix} F_i(P)\\ M_i(P)\end{pmatrix}$ of the six limbs is equal to five, then there is a line of twists reciprocal to all $\begin{pmatrix} F_i(P)\\ M_i(P)\end{pmatrix}$.

Let $A=(a_{i,j})_{i,j=1,\ldots,n}$ be a square matrix of size $n$ whose entries are polynomials in $X=(x_1,\ldots,x_d)$ with coefficients in the field $k$.
\begin{lemma}\label{corang1}
If $\det(A)=0$ is a nonsingular hypersurface in $d$-dimensional affine space $k^d$, then the matrix $A$ has rank $n-1$ at every point of this hypersurface. Moreover, if $\mathrm{sol}(A)$ denotes the line of solutions $S$ of $AS=0$ for every $A$ such that $\det(A)=0$, then $A\mapsto \mathrm{sol}(A)$ is a regular mapping from the hypersurface $\det(A)=0$ to the projective space $\mathbb{P}^{n-1}(k)$.
\end{lemma}
\noindent \textbf{Proof.}
The derivative of $\det(A)$ with respect to $x_\ell$ can be expressed as
\begin{equation}
\frac{\partial \det(A)}{\partial x_\ell}=\sum_{j=1}^n\sum_{i=1}^n \mathrm{cof}(A)_{i,j} \frac{\partial a_{i,j}}{\partial x_\ell}\;,
\end{equation}
where $\mathrm{cof}(A)_{i,j}$ is the cofactor of $a_{i,j}$ in matrix $A$. Assume that there is a point in the affine space of dimension $d$ where the rank of $A$ is $<n-1$. Then all minors of dimension $n-1$ of the matrix $A$ vanish at this point and the formula above shows that all partial derivatives of $\det(A)$ also vanish at this point. This contradicts the assumption that $\det(A)=0$ is a nonsingular hypersurface.\par
Let $U_i$ be the Zariski open subset of the nonsingular hypersurface $\det(A)=0$ where all rows of $A$ but the $i$-th are linearly independent. By the first part of the lemma, $(U_i)_{i=1,\ldots,n}$ is a Zariski open cover of $\det(A)=0$. For any $A$ in $U_i$, the line of solutions of $AS=0$ is spanned by the vector with coordinates $\mathrm{cof}(A)_{i,j}$ for $j=1,\ldots,n$; these cofactors are polynomial in the coefficients of $A$, hence also polynomials in $X$. In the end, we obtain a regular mapping from the hypersurface $\det(A)=0$ in $k^d$ to $\mathbb{P}^{n-1}(k)$. \quad $\square$

The aim of this paper is to prove the rationality of $\Sing$ for a general Gough-Stewart platform, without assuming planarity of the base or of the platform.
\begin{proposition} Assume that $\Sigma$ is a nonsingular hypersurface of $\R(\SO(3))^3$. Then, for every singular configuration in $\Sigma$, the system of Pl\"ucker coordinates of the six limbs is of rank $5$ in the space of screws. The mapping which associates, to each singular configuration, the line of twists reciprocal to the Pl\"ucker coordinates of the six limbs is a regular mapping from $\Sigma$ to $\mathbb{P}^5(\R(\SO(3)))$, which we denote by $\mathrm{Rec}$.
\end{proposition}

\noindent \textit{Proof.}
Lemma \ref{corang1} shows that the rank of the matrix $\Jac$ of the Pl\"ucker coordinates of the six limbs is indeed $5$ for each singular configuration $P$. The matrix $\Jac$ is the matrix of the system expressing the reciprocity of the twist $\begin{pmatrix}
\Omega\\ V
\end{pmatrix}$ to the systems of Pl\"ucker coordinates $\begin{pmatrix}
F_i(P)\\ M_i(P)
\end{pmatrix}$ of the six limbs. The second part of Lemma \ref{corang1} shows that the mapping $\mathrm{Rec}$ associating to each singular configuration $P$ the line of screws reciprocal the the Pl\"ucker coordinates of the six limbs is indeed a regular mapping from $\Sigma$ to $\mathbb{P}^5{\R(\SO(3))}$.\quad $\square$\par\medskip

We are going to show in the following that the image of $\mathrm{Rec}$ is contained in a quadric surface in a $3$-dimensional subspace of $\mathbb{P}^5(\R(\SO(3)))$ and that $\mathrm{Rec}$ is a birational equivalence with this quadric surface.

\subsection{The image of $\mathrm{Rec}$ is contained in a quadric surface}

The equations expressing that the twist $\begin{pmatrix}
\Omega\\ V
\end{pmatrix}$ is reciprocal to the systems of Pl\"ucker coordinates of the limbs are
\begin{equation}
(C_i+P)\cdot V + (A_i\times (C_i+P))\cdot \Omega=0,\quad i=1,\ldots,6.
\end{equation}
This is equivalent to:
\begin{equation}\label{Eqreci1}
V\cdot C_i+[\Omega ,A_i,C_i]+(\Omega \times A_i+V)\cdot P=0,\quad i=1,\ldots,6,
\end{equation}
where $[\Omega ,A_i,C_i]$ is the mixed product of vectors $\Omega ,A_i,C_i$. 

We are going to eliminate $P$ from this system of equations. The first step is to obtain two linear equations in $V,\Omega$ where $P$ does not appear. In equation $(\ref{Eqreci1})_1$, $A_1=C_1=0$ so $(\ref{Eqreci1})_1$ is $V\cdot P=0.$ Subtracting this first equation from the other ones we get:
\begin{equation}\label{Eqreci2}
V\cdot C_i+[\Omega ,A_i,C_i]+(\Omega \times A_i)\cdot P=0,\quad i=2,\ldots,6.
\end{equation}
The family of vectors $A_i, i=2,\ldots,6$ has rank at most $3.$ Up to reordering the indices, we may assume that $A_5$ and $A_6$ are linear combinations of $A_2, A_3$ and $A_4.$ Set $A_5=\alpha_2 A_2+\alpha_3 A_3+ \alpha_4 A_4$ and $A_6=\beta_2 A_2+\beta_3 A_3+ \beta_4 A_4.$  Computing $(\ref{Eqreci2})_5-\sum_{i=2}^{4}\alpha_i (\ref{Eqreci2})_i$ and $(\ref{Eqreci2})_6-\sum_{i=2}^{4}\beta_i (\ref{Eqreci2})_i$  we obtain respectively
\begin{equation}\label{eqlin1}
V\cdot (C_5-\sum_{i=2}^{4}\alpha_iC_i)+[\Omega ,A_5,C_5]
-\sum_{i=2}^{4}\alpha_i [\Omega ,A_i,C_i]=0
\end{equation}
and
\begin{equation}\label{eqlin2}
V\cdot (C_6-\sum_{i=2}^{4}\beta_iC_i)+[\Omega ,A_6,C_6]
-\sum_{i=2}^{4}\beta_i [\Omega ,A_i,C_i]=0
\end{equation}
which are linear homogeneous equations in $\Omega , V.$

The second step in the elimination of $P$ is to obtain a quadratic equation in $V,\Omega$ without $P$.
Set $\ell_2=[\Omega ,A_3,A_4],\, \ell_3=[\Omega ,A_4,A_2],\, \ell_4=[\Omega ,A_2,A_3].$
\begin{lemma}
$\ell_2(\Omega \times A_2)+ \ell_3(\Omega \times A_3)+\ell_4(\Omega \times A_4)=0.$
\end{lemma}
\noindent \textit{Proof.} It suffices to check this formula of vector algebra assuming $A_2,A_3,A_4$ linearly independent. The dot product of the left handside with any of $A_2,A_3,A_4$ is easily seen to be zero, hence, it is the zero vector.\quad$\square$\par\medskip

The linear combination of  $(\ref{Eqreci2})_2,(\ref{Eqreci2})_3,(\ref{Eqreci2})_4$ with coefficients $\ell_2,\ell_3,\ell_4$ yields the following homogeneous quadratic equation in $\Omega ,V$:
\begin{multline}\label{eqquad}
[\Omega ,A_3,A_4](V\cdot C_2+[\Omega ,A_2,C_2])+[\Omega ,A_4,A_2](V\cdot C_3+[\Omega ,A_3,C_3])\\
+[\Omega ,A_2,A_3](V\cdot C_4+[\Omega ,A_4,C_4])=0\;.
\end{multline}
Equations (\ref{eqlin1}), (\ref{eqlin2}) and (\ref{eqquad}) are the homogeneous equations of a quadric surface $Q$ in $\mathbb{P}^5(\R(\SO(3)))$, the projectivisation of the space of twists. We shall check in the case study that this projective quadric surface is nonsingular for a generic Gough-Stewart platform. We have proved

\begin{proposition} The image of the regular mapping $\mathrm{Rec}$ is contained in the quadric surface $Q$ given by equations (\ref{eqlin1}), (\ref{eqlin2}) and (\ref{eqquad}).
\end{proposition}

\subsection{The birational equivalence}

We obtained the system of equations (\ref{eqlin1}), (\ref{eqlin2}) and (\ref{eqquad}) from the system of equations in $P$, $V$ and $\Omega $ expressing that the screw $\begin{pmatrix}
\Omega \\ V
\end{pmatrix}$ is reciprocal to the Pl\"ucker coordinates of the six limbs for the singular configuration $P$. Let $\Omega ,V$ be the homogeneous coordinates of a point on the quadric surface $Q$ such that $[V,\Omega \times A_2, \Omega \times A_3]\neq 0$. Generically, the quadric $Q$ is nonsingular and the set of $\Omega ,V$ such that $[V,\Omega \times A_2, \Omega \times A_3]\neq 0$ is a dense Zariski open subset of $Q$. The condition $[V,\Omega \times A_2, \Omega \times A_3]\neq 0$ implies that $\ell_4=[\Omega ,A_2,A_3]\neq 0$ ; indeed $[\Omega ,A_2,A_3]= 0$ implies that $\Omega \times A_2$ and $\Omega \times A_3$ are colinear.
Since $\ell_4\neq 0$, the system of equations $(\ref{Eqreci1}_i)$ for $i=1,\ldots,6$ is equivalent to the system of equations $(\ref{Eqreci1}_1, \ref{Eqreci2}_2, \ref{Eqreci2}_3,\ref{eqlin1},\ref{eqlin2},\ref{eqquad})$. The variables $P$ appear only in the first three equations of this system and these equation can be rewritten as
\begin{equation}\label{systP}
\begin{aligned}
V\cdot P&=0\\
(\Omega \times A_2)\cdot P&=-V\cdot C_2-[\Omega ,A_2,C_2]\\
(\Omega \times A_3)\cdot P&=-V\cdot C_3-[\Omega ,A_3,C_3]\;.
\end{aligned}
\end{equation}
This is a linear system in $P$ which is a Cramer system since we have taken $(\Omega ,V)$ in the quadric surface $Q$ such that its determinant $[V,\Omega \times A_2, \Omega \times A_3]$ does not vanish. Hence there is a unique $P$ satisfying the system (\ref{systP}), and this $P$ is given by the Cramer formulas, i.e. as a rational function (of degree 3) in $\Omega ,V$. This $P$ is on the surface $\Sigma$, since $P,V,\Omega $ satisfy the reciprocity equations, which implies that the rank of the system of Pl\"ucker coordinates of the limbs for $P$ is $<6$. So we get a rational mapping, which we denote by $\mathrm{Pos}$, from $Q$ to $\Sigma$. The uniqueness of the solution for the Cramer system (\ref{systP}) implies that the composition $\mathrm{Pos}\circ \mathrm{Rec}$ is the identity of $\Sigma$. The fact that there is only a line of screws reciprocal to the Pl\"ucker coordinates of the limbs for a configuration $P$ in $\Sigma$ shows that the composition $\mathrm{Rec}\circ \mathrm{Pos}$ is the identity on $Q$.\par\medskip

In conclusion, we have proved the following theorem.

\begin{theorem}\label{birat}
The regular mapping $\mathrm{Rec}: \Sigma \to Q$  is a birational equivalence to the quadric surface $Q\in \mathbb{P}^5(\R(\SO(3)))$, defined over $\R(\SO(3))$.
\end{theorem}

\begin{corollary}\label{sigmarat}
The cubic surface $\Sigma$ is rational over $\R(\SO(3))$.
\end{corollary}

\noindent \textit{Proof.} The cubic surface $\Sigma$ has the $\R(\SO(3))$-rational point $P=0$. Hence, the quadric surface $Q$ has the $\R(\SO(3))$-rational point $\mathrm{Rec}(0)$. It follows that $Q$ is rational over $\R(\SO(3))$ and since $\Sigma$ is $\R(\SO(3))$-birational to $Q$, it is also rational over $\R(\SO(3))$.\quad$\square$

\begin{theorem}
For a generic Gough-Stewart platform, the hypersurface $\Sing\subset \SE(3)$ of singular configurations is rational over $\R$.
\end{theorem}

\noindent \textit{Proof.} As already mentioned, this is a consequence of Corollary \ref{sigmarat} and of the fact that $\SO(3)$ is rational over $\R$.\quad$\square$

\subsection{Case study continued}\label{Case2}
We continue here the computations for the cubic surface $\Sing_R$ where the geometry of the Gough-Stewart platform and the orientation are those of section \ref{Case}.
Let us denote by $(\omega_1,\omega_2,\omega_3,v_1,v_2,v_3)$ the coordinates of a twist. We compute the variety of reciprocal twists as follows: we form the ideal $\mathfrak{Rec}$ generated by the reciprocal products of the twist with the Pl\^ucker coordinates of the limbs and the Jacobian determinant, then we eliminate the variables $x,y,z$. The ideal obtained is the homogeneous ideal generated by
\begin{equation}
\begin{aligned}
&72\,\omega_1-42\,\omega_2+16\,v_1+21\,v_3+35\,v_2,\\
&18\,\omega_3+30\,\omega_2-26\,v_1-15\,v_3+5\,v_2,\\
&168\,\omega_2^2+80\,\omega_2\,v_1+522\,\omega_2\,v_3-296\,v_1\,v_3-159\,v_3^2\\
&\qquad{}-518\,\omega_2\,v_2+264\,v_1\,v_2+212\,v_3\,v_2-165\,v_2^2\\
\end{aligned}
\end{equation}
There are two linear equations and one quadratic: this is the ideal of a quadric surface in $\Proj^5(\R)$.

We compute formulas for the reciprocal twist (only defined up to a scalar factor, as a point in $\Proj^5(\R)$) in the following way: eliminating variables $\omega_1,\omega_2,\omega_3,v_3$ from the ideal $\mathfrak{Rec}$, we obtain an ideal which contains a polynomial which is linear homogeneous in $v_1,v_2$ and quadratic in $x,y,z$. This gives $v_2/v_1$ as a rational function of degree $2$ in $x,y,z$. Proceeding in the same way to compute $v_3/v_1,\omega_1/v_1,\omega_2/v_1,\omega_3/v_1$ and chasing denominator (the same for all expressions), we arrive at
\begin{equation}\label{reciscrew}
\begin{aligned}
v_1&=-80\,x^2+55\,x\,y-109\,x\,z+140\,y^2-14\,y\,z-84\,z^2\\
&\qquad{}+376\,x+774\,y+20\,z+1392\\
v_2&= 52\,x^2-131\,x\,y+164\,x\,z+68\,y^2-62\,y\,z+24\,z^2\\
&\qquad{}-473\,x+78\,y-130\,z-426\\
v_3&= 156\,x^2-245\,x\,y+180\,x\,z-36\,y^2-102\,y\,z+24\,z^2\\
&\qquad{}-663\,x-578\,y-234\,z-1410\\
\omega_1&= -53\,x^2+92\,x\,y-108\,x\,z+338\,x+248\,y-48\,z+736\\
\omega_2&= -53\,x\,y+92\,y^2-108\,y\,z-3\,x+496\,y-300\,z+732\\
\omega_3&= -53\,x\,z+92\,y\,z-108\,z^2+127\,x-212\,y+370\,z-266\\
\end{aligned}
\end{equation}
These formulas are degree $2$ polynomials in $x,y,z$. The formulas define a point in $\Proj^5(\R)$ except when they all vanish. We compute the ideal generated by these quadratic polynomials in $x,y,z$ and the cubic equation of the surface of singularities. We find that this ideal is precisely the ideal of the two lines corresponding to factor $F_2$.

Nevertheless, the regular mapping which associates to a point $(x,y,z)$ of the cubic singularity surface, a point of $\Proj^5(\R)$ representing the line of reciprocal twists is also well defined on these two lines. Indeed, at every point in the surface of singularities there is a nonzero cofactor in the Jacobian matrix, and the cofactors on the same row as this nonzero cofactor are the coordinates of a nonzero reciprocal twist. These cofactors are polynomials of degree $3$ in $x,y,z$. We shall explain in section \ref{rectwistinf} why we obtain degree $2$ in the formulas above instead of degree $3$. 

\begin{figure}[h]\label{Quadricandpoints}
\caption{The quadric surface of reciprocal twists in blue}
\begin{center}
\includegraphics[scale=.5]{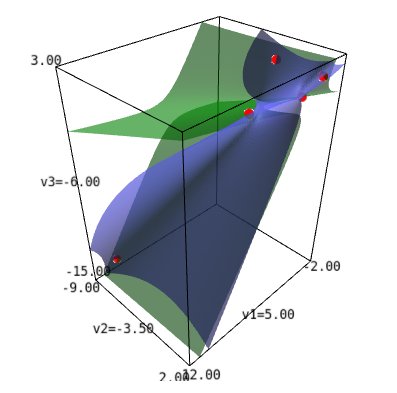}
\end{center}
\end{figure}

We next compute the image of the five lines corresponding to factor $F_5$ by the reciprocal twist mapping. This is again the computation of an elimination ideal, and the computed ideal is zero-dimensional, of degree $5$. It is actually the ideal of five distinct points of the quadric, all real in the case under consideration. The reciprocal twists corresponding to these five points are self-reciprocal. The figure \ref{Quadricandpoints} is in the $3$-space containing the quadric, in its affine chart given by $\omega_2=1$, using coordinates $v_1,v_2,v_3$. The quadric image of the reciprocal twist mapping is in blue, the five points in red, and the intersection of the quadric of self-reciprocal twists with the $3$-space is in green.

It is also possible to compute formulas for the inverse rational mapping from the quadric surface in $\Proj^5(\R)$ to $\Sing_R$: for instance, eliminating the variables $y,z$ in the ideal $\mathfrak{Rec}$, we obtain $x$ as a rational function of degree $2$ in the $V,\Omega$. Actually we compute formulas for the rational mapping to the projective closure $(\Sing_R)^h$ in $\Proj^3(\R)$. We use homogeneous coordinates $(w:x:y:z)$ for $\Proj^3(\R)$, with $w=0$ as plane at infinity.
\begin{equation}\label{invrec}\begin{aligned}
w = {}&60v_3^2+83v_1v_3+42\omega_2v_1-20v_3v_2-35v_1v_2
+72\omega_2v_2-120\omega_2v_3-16v_1^2
\\
x = {}& 162\omega_2v_3+54\omega_2v_2-48v_1v_3-48v_1v_2-213v_2^2
+63v_3^2-42v_3v_2\\
y = {}& -54\omega_2v_1+48v_1^2+288\omega_2v_3-357v_1v_3-144v_3^2
+213v_1v_2-48v_3v_2\\
z = {}& -162\omega_2v_1-288\omega_2v_2-63v_1v_3+399v_1v_2
+48v_2^2+48v_1^2+144v_3v_2
\end{aligned}
\end{equation}
Note that these formulas are of degree $2$, whereas the Cramer formulas invoked in the argument for birationality would give degree $3$. It is possible to check that formulas (\ref{reciscrew}) and (\ref{invrec}) give actually a birational isomorphism between $(\Sing_R)^h$ and the quadric surface in $\Proj^5(\R)$. One can also compute the indetermination points of formulas (\ref{invrec}) on the quadric surface: one finds five points, which are precisely the image of the five lines corresponding to $F_5$. This indicates that the ``reciprocal twist" regular mapping from $\Sigma$ to $Q$ may be extended to a regular mapping from the projective closure $\Sigma^h$ to $Q$, which is the blowing-up of $Q$ in five points. We shall prove this in the following section after investigating the points at infinity of the singularity locus, i.e. $\Sigma^h\setminus \Sigma$.  

\section{The projective closure of the singularity surface}\label{Sec4}

We already know that the projective closure $\Sigma^h$ of $\Sigma$ in $\Proj^3(\R(\SO(3)))$ is a nonsingular projective cubic surface for a generic Gough-Stewart platform. Indeed, we checked in section \ref{Case} that $(\Sing_R)^h$ is nonsingular. 

\subsection{Singularity condition at infinity\label{Secsinginf}}
 
The formula (\ref{jacobien}) can be rewritten as a sum indexed by the symmetric group $\mathfrak S_6$, with $\epsilon(\sigma)$ denoting the signature of a permutation $\sigma$:

\begin{equation}\label{jacobienperm}
\begin{aligned}
\det(\Jac) = &\frac1{36}\;\sum_{\sigma \in \mathfrak S_6} \epsilon(\sigma)\, [C_{\sigma(1)}+P, C_{\sigma(2)}+P, C_{\sigma(3)}+P]\\ &\qquad [A_{\sigma(4)}\times (C_{\sigma(4)}+P), A_{\sigma(5)}\times (C_{\sigma(5)}+P), A_{\sigma(6)}\times (C_{\sigma(6)}+P)]\;.
\end{aligned}
\end{equation}

The equation of the projective closure $\Sigma^h$ is the homogeneization of $\det(\Jac)$, say with homogenization variable $w$. The equation of its intersection with the plane at infinity $w=0$ is the homogeneous part of degree $3$ of $\det(\Jac)$. Using multilinearity in (\ref{jacobienperm}), we find that this homogeneous part is

\begin{equation}
H_3 = \frac1{4}\;\sum_{\sigma \in \mathfrak S_6} \epsilon(\sigma)\, [C_{\sigma(1)}, C_{\sigma(2)},P]\,[A_{\sigma(4)}\times C_{\sigma(4)}, A_{\sigma(5)}\times P, A_{\sigma(6)}\times P]\;.
\end{equation}
Since
\begin{equation}
[A_{\sigma(4)}\times C_{\sigma(4)}, A_{\sigma(5)}\times P, A_{\sigma(6)}\times P]=[A_{\sigma(4)},  C_{\sigma(4)}, P] [A_{\sigma(5)},  A_{\sigma(6)}, P]\;, 
\end{equation}
we arrive to
\begin{equation}
\begin{aligned}
H_3 = {}&\frac1{4}\;\sum_{\sigma \in \mathfrak S_6} \epsilon(\sigma)\, [C_{\sigma(1)}, C_{\sigma(2)},P]\,[A_{\sigma(4)},  C_{\sigma(4)}, P]\, [A_{\sigma(5)},  A_{\sigma(6)}, P]\\
= {}&\frac{-1}{4} \sum_{\sigma\in \mathfrak{S}_6} \varepsilon(\sigma) [A_{\sigma(1)}, A_{\sigma(2)}, P]\,[A_{\sigma(3)}, C_{\sigma(3)}, P]\,[C_{\sigma(4)}, C_{\sigma(5)}, P]
\end{aligned}
\end{equation}

In this formula we can replace the $C$ with the $B$:

\begin{equation}
H_3=\frac{-1}{4} \sum_{\sigma\in \mathfrak{S}_6} \varepsilon(\sigma) [A_{\sigma(1)}, A_{\sigma(2)}, P][A_{\sigma(3)}, B_{\sigma(3)}, P][B_{\sigma(4)}, B_{\sigma(5)}, P]
\end{equation}

The equality is proved using $C_i=B_i-A_i$ and Grassmann-Pl\"ucker relations such as:
\begin{equation}
\begin{aligned}
{}[A_{\sigma(1)}, A_{\sigma(2)}, P][A_{\sigma(4)}, A_{\sigma(5)}, P] - [A_{\sigma(1)}, A_{\sigma(4)}, P][A_{\sigma(2)}, A_{\sigma(5)}, P]&{}\\{} +[A_{\sigma(1)}, A_{\sigma(5)}, P][A_{\sigma(2)}, A_{\sigma(4)}, P]&=0
\end{aligned}
\end{equation}

We summarize what we have seen in the following

\begin{proposition}\label{Singinf}
The singularities at infinity of the Gough-Stewart platform in the orientation $R$ are the directions of $P=\T{\begin{pmatrix}
x&y&z
\end{pmatrix}}$
satisfying
$$\sum_{\sigma\in \mathfrak{S}_6} \varepsilon(\sigma) [A_{\sigma(1)}, A_{\sigma(2)}, P][A_{\sigma(3)}, B_{\sigma(3)}, P][B_{\sigma(4)}, B_{\sigma(5)}, P]=0\;,$$
where $B_i=Rb_i$.
\end{proposition}

We can reformulate this condition as a condition concerning two planar hexagons. Let $\pi_P$ the orthogonal projection in the direction of $ P. $ Then for $i=1,\ldots,6,$ we denote $\pi_P(A_i)=\alpha_i$ the projection of the vertices of the base and $\pi_P(B_i)=\beta_i$ the projection of the vertices of the platform after rotation $R$. So we have two hexagons in the plane orthogonal to $P.$ The condition of singularity at infinity is given by:
\begin{equation}\label{hexagons}
\sum_{\sigma\in \mathfrak{S}_6} \varepsilon(\sigma) [\alpha_{\sigma(1)}, \alpha_{\sigma(2)}][\beta_{\sigma(3)}, \beta_{\sigma(4)}][\alpha_{\sigma(5)}, \beta_{\sigma(5)}]=0\;,
\end{equation}
where the bracket denotes here $2\times 2$ determinants. This condition is invariant under the following transformations: 
\begin{enumerate}
\item The same affine transformation applied to all $\alpha_i$'s and $\beta_i$'s. The determinants are multiplied by the determinant of linear part of the transformation.
\item  A homothety (of ratio $k\not= 0$ ) on one of the hexagons. The product of determinants are multiplied by $k^3.$
\item A translation (of vector $u$) on one of the hexagons. The invariance follows from Pl\"ucker relations.
\item Replace the $\beta_i'$ s by any linear combination of $\beta_i$ and $\alpha_i.$  
\end{enumerate}
In addition, when the hexagons have three distinct vertices in common ($\alpha_1=\beta_1, \alpha_2=\beta_2,\alpha_3=\beta_3$), the condition is satisfied if and only if the three lines $(\alpha_i,\beta_i)$ for $i=4,5,6$ are concurrent or parallel.

\subsection{Reciprocal twist at infinity}\label{rectwistinf}

Recall that the twist $\begin{pmatrix}
\Omega \\ V
\end{pmatrix}$ is reciprocal to the system of Pl\"ucker coordinates of the six limbs if and only if
\begin{equation}
(C_i+P)\cdot V + (A_i\times (C_i+P)) \cdot \Omega =0,\quad i=1,\ldots,6.
\end{equation}
The matrix of this system is  $\Jac$.
In a singular configuration, the rank of this system is generically $5.$
In this case, if the rows of $\Jac$ with index $\neq i$ are linearly independant, the line of solutions of the system of equations is spanned by the vector whose coordinates are $\mathrm{cof}(\Jac)_{i,j}$ for $j=1,\ldots,6$. This gives formulas of degree $3$ in $P$. The line of solutions is also spanned by the vector whose coordinates are $\sum_{i=1}^6\mathrm{cof}(\Jac)_{i,j}$ for $j=1,\ldots,6$, if this vector is nonzero. It turns out that these sums of cofactors are polynomials of degree $2$ in $P$; this explains what we have seen in the case study in formulas (\ref{reciscrew}).

The sum by columns of the cofactors in the first three columns is the vector

\begin{equation}\begin{aligned}
T_1 &=  \sum_{1\le i_1< i_2< i_3\le 6} (-1)^{i_1+ i_2+ i_3} \left(\sum_{\mathrm{cyc}} (C_{i_2}+P)\times (C_{i_3}+P)\right)\\
&\qquad\qquad [(A_{j_1}\times (C_{j_1}+P), A_{j_2}\times (C_{j_2}+P), A_{j_3}\times (C_{j_3}+P)]
\end{aligned}\end{equation}
where $j_1<j_2<j_3$ are the integers between $1$ and $6$ different from $i_1,i_2,i_3$ and the cyclic sum is taken over the powers of the cycle $(i_1,i_2,i_3)$. It will be more convenient to write the sum $T_1$ as a sum over all permutations in $\mathfrak S_6$.

\begin{equation}
\begin{aligned}
T_1&=\frac{1}{12}\sum_{\sigma \in \mathfrak S_{6}}\Bigg (\varepsilon(\sigma)\, (C_{\sigma(2)}+P) \times (C_{\sigma(3)}+P) \\
&\qquad [A_{\sigma(4)}\times ( C_{\sigma(4)}+P), A_{\sigma(5)}\times ( C_{\sigma(5)}+P), A_{\sigma(6)}\times ( C_{\sigma(6)}+P) ] \Bigg)
\end{aligned}
\end{equation}
We develop this expression using multilinearity in the cross-product and in the mixed product. We remark that all terms containing only four among the six indices $\sigma(i)$ disappear in the sum over all permutations (due to the transposition on the missing indices). We remark also that the mixed product
$$[A_{\sigma(4)}\times P, A_{\sigma(5)}\times P, A_{\sigma(6)}\times P ]
$$
is zero because the three vectors are linearly dependant. This explains why all terms of degree 3 disappear in the sum and why the only remaining terms of degree 2 are 
\begin{equation}\label{torsreciinf1}\begin{aligned}
V_{\infty}&=\frac{1}{4}\sum_{\sigma \in \mathfrak S_6}\varepsilon(\sigma)\; [A_{\sigma(4)}\times  C_{\sigma(4)}, A_{\sigma(5)}\times P, A_{\sigma(6)}\times P ]\; C_{\sigma(2)} \times C_{\sigma(3)}\\
&=\frac{1}{4}\sum_{\sigma \in \mathfrak S_{6}}\varepsilon(\sigma)
\;[A_{\sigma(4)},C_{\sigma(4)},P]\, [A_{\sigma(5)},A_{\sigma(6)},P]\; C_{\sigma(2)} \times C_{\sigma(3)}\\
&=\frac{-1}{4}\sum_{\sigma \in \mathfrak S_{6}}\varepsilon(\sigma)\;  [A_{\sigma(1)}, A_{\sigma(2)}, P ] \, [A_{\sigma(3)}, C_{\sigma(3)}, P ]\; C_{\sigma(4)} \times C_{\sigma(5)}
\end{aligned}
\end{equation}

The sum by columns of the cofactors in the last three columns is the vector

\begin{equation}
\begin{aligned}
T_2&=\sum_{1\le i_1< i_2< i_3\le 6} (-1)^{i_1+ i_2+ i_3} \bigg( [C_{i_1}+P, C_{i_2}+P, C_{i_3}+P]\\
&\qquad\qquad\qquad \sum_{\mathrm{cyc}} (A_{j_1}\times ( C_{j_1}+P))\times (A_{j_2}\times ( C_{j_2)}+P)\bigg)
\end{aligned}
\end{equation}
where $j_1<j_2<j_3$ are the integers between $1$ and $6$ different from $i_1,i_2,i_3$ and the cyclic sum is taken over the powers of the cycle $(j_1,j_2,j_3)$. Here also it will be more convenient to write the sum $T_2$ as a sum over all permutations in $\mathfrak{S}_6$.

\begin{equation}
\begin{aligned}
T_2&=\frac1{12}\sum_{\sigma \in \mathfrak S_{6}}\varepsilon(\sigma)[C_{\sigma(1)}+P, C_{\sigma(2)}+P, C_{\sigma(3)}+P]\\
&\qquad\qquad\qquad (A_{\sigma(4)}\times ( C_{\sigma(4)}+P))\times (A_{\sigma(5)}\times ( C_{\sigma(5)}+P)
\end{aligned}
\end{equation}

We develop this expression using multilinearity in the mixed product and in the cross-product. Here also, terms where only
four indices among the $\sigma(i)$ are present disappear in the sum over all permutations. This explains why all terms of degree $3$ in $P$ disappear, and why the terms of degree $2$ in $P$ which remain are only:

 \begin{equation}\label{torsreciinf2}
 \begin{aligned}
\Omega _{\infty}&=\frac1{12}\Bigg(\sum_{\sigma \in \mathfrak S_{6}}\varepsilon(\sigma) [C_{\sigma(1)}, C_{\sigma(2)}, C_{\sigma(3)}]\,(A_{\sigma(4)}\times P)\times (A_{\sigma(5)}\times P)\Bigg)\\
&=\frac1{12}\left(\sum_{\sigma \in S_{6}}\varepsilon(\sigma) [C_{\sigma(1)}, C_{\sigma(2)}, C_{\sigma(3)}] [A_{\sigma(4)}, A_{\sigma(5)}, P ]\right) P 
\end{aligned}
\end{equation}
Note that $\Omega_{\infty}$ is the product of the vector $P$ by a linear form in $P$; this confirms the result of the computation in case study \ref{Case2}.

The notations $V_\infty$ and $\Omega_\infty$ are explained by the fact that $\begin{pmatrix} \Omega_\infty\\ V_\infty\end{pmatrix}$ is actually the reciprocal twist at infinity, as explained in the following

\begin{proposition}\label{Rech}
For a generic Gough-Stewart platform, the reciprocal twist mapping $\mathrm{Rec}: \Sigma\to Q$ can be extended to a regular mapping $\mathrm{Rec}^h: \Sigma^h\to Q$. The image by $\mathrm{Rec}^h$ of a point at infinity in the direction $P$ in $\Sigma^h$ is given by the homogeneous second degree formulas (\ref{torsreciinf1}) for $V_{\infty}$ and (\ref{torsreciinf2}) for $\Omega_{\infty}$, provided that $V_{\infty}$ and $\Omega_{\infty}$ do not both vanish at $P$.
\end{proposition}

\noindent \textit{Proof.} We know that $\begin{pmatrix}
T_2\\ T_1
\end{pmatrix}$ give degree $2$ formulas for $\mathrm{Rec}$ at points of $\Sigma$ where $T_1$ and $T_2$ do not both vanish. Hence, the homogeneous quadratic part $\begin{pmatrix}
\Omega_{\infty}\\ V_{\infty}
\end{pmatrix}$ give a regular extension of $\mathrm{Rec}$ at points at infinity of $\Sigma^h$ in the direction of $P$, when $\Omega_\infty$ and $V_\infty$ do not both vanish.
We note that, in the plane at infinity, $\Omega_\infty$ vanish on a line (whose equation is the linear form $L(P)$ such that $\Omega_\infty=L(P)\,P$), while each component of $V_\infty$ vanishes on a conic. We check in the case study \ref{Case3} that these curves have two points in common which are on the cubic curve of the points at infinity of $\Sigma^h$; these two points are, of course, the points at infinity of the two lines on $\Sing_R^h$ corresponding to $F_2$. We check also that the formulas for the reciprocal twist at infinity given by the homogeneous part of degree $3$ of the cofactors of $\Jac$ on the first line do not vanish at these two points. This shows that for a generic Gough-Stewart platform, the regular extension $\mathrm{Rec}^h$ is well defined at every point at infinity of $\Sigma^h$: when both  $\Omega_\infty$ and $V_\infty$ vanish, it is given by the homogeneous part of degree $3$ of the cofactors of $\Jac$ on the first line.\quad$\square$\par\medskip

\begin{proposition}\label{selfreci}
For a generic Gough-Stewart platform, the reciprocal twists at points at infinity of $\Sigma^h$ in the direction $P$ are all self-reciprocal. They are the twists of a rotation with axis parallel to $P$, or of a translation orthogonal to $P$. 
\end{proposition}

\noindent \textit{Proof.} From formulas (\ref{torsreciinf1}) for $V_{\infty}$ and Proposition (\ref{Singinf}), we obtain $V_\infty\cdot P=0$ for all point at infinity of $\Sigma^h$ in direction $P$. Since $\Omega_\infty=L(P) P$, the statements of the Proposition follow. \quad$\square$\par\medskip

\subsection{$\Sigma^h$ as the blowing-up in five points of a quadric surface}

\begin{theorem}
For a generic Gough-Stewart platform, the reciprocal twist regular mapping $\mathrm{Rec}^h: \Sigma^h\to Q$ is the blowing-up of the projective quadric surface $Q$ in five points. The five exceptional divisors are five lines on $\Sigma^h$, forming a $S_5$ intersecting both lines of a $S_2$ (using notations of Theorem \ref{SwDy}), both $S_5$ and $S_2$ being defined over $\R(\SO(3))$. A point of $\Sigma$ belongs to one of the lines of the $S_5$ if and only if its reciprocal twist is the twist of a pure rotation; the axis of this rotation is then a line intersecting the six limbs of the platform.
\end{theorem}

\noindent \textit{Proof.} 
We know from Proposition \ref{Rech} that $\mathrm{Rec}^h : \Sigma^h\to Q$ is a regular on the whole of $\Sigma^h$, and we know from Theorem \ref{birat} that it is a birational equivalence. From this we deduce that $\mathrm{Rec}^h$ is the composition of a sequence of blowing-ups at points (see for instance \cite{Har}, corollary 5.4 p.411). Since the divisor class group of a nonsingular projective cubic surface is $\mathbb Z^7$ (\cite{Har} Proposition 4.8 p.401) and the divisor class group of a nonsingular projective quadric surface is $\mathbb Z^2$, $\mathrm{Rec}^h$ is the blowing-up of $Q$ in five points. These five points are distinct points of $Q$, as checked in the case study, so the five exceptional divisors have self-intersection equal to $-1$ and are mutually skew lines on $\Sigma^h$. Since $\mathrm{Rec}^h$ is defined over $\R(\SO(3))$, the five exceptional divisors form a $S_5$ over $\R(\SO(3))$. There are two possible types of $S_5$, one of them characterized as the five lines intersecting both lines of an $S_2$ (see \cite{SD}). Actually, the five exceptional divisors obtained by blowing up five points in a quadric form a $S_5$ of this type; we defer the proof of this fact to Lemma \ref{eclatement} below.

Since the reciprocal twist is constant along each line of the $S_5$, it is equal to the reciprocal twist of the point at infinity of this line. Hence, by  Proposition \ref{selfreci}, it is the twist of a rotation or of a translation. If the reciprocal twist at $P\in \Sigma$ were the twist $\begin{pmatrix} 0\\V\end{pmatrix}$ of a translation in the direction of vector $V$ then, for every vector $U$ orthogonal to $V$, the twist $\begin{pmatrix} 0\\V\end{pmatrix}$ would be reciprocal to the systems of Pl\"ucker coordinates of all limbs in position $P+U$ and hence $P+U$ would be in $\Sigma$. This cannot be since $\Sigma$ contains no plane. This shows that the reciprocal twist at each point of $\Sigma$ belonging to a line of the $S_5$ is the twist of a pure rotation; the reciprocity conditions means that every limb of the platform intersects (or is parallel to) the axis of this rotation.

Reciprocally, suppose that the reciprocal twist at a point $P\in \Sigma$ is the twist  $\begin{pmatrix} \Omega\\V\end{pmatrix}$ of a pure rotation, i.e. $\Omega\neq 0$ and $\Omega\cdot V=0$. Then, for all scalars $\lambda$, $P+\lambda \Omega$ is still in $\Sigma$, with reciprocal twist $\begin{pmatrix} \Omega\\V\end{pmatrix}$. Hence the line of $P+\lambda \Omega$ is one of the exceptional divisors in the $S_5$. \quad$\square$\par\medskip

\begin{lemma}\label{eclatement} The five exceptional divisors obtained by blowing up a quadric surface in five distinct points form an $S_5$ with two transversals forming an $S_2$.
\end{lemma}

\noindent \textit{Proof.} 
Let $p_1,\ldots,p_4,q$ be the five distinct points on the quadric surface. Let $d_5,d_6$ be the two lines on the quadric through $q$. Blow up $q$; the strict transforms of $d_5$ and $d_6$ (which we still denote by $d_5$ and $d_6$) have self-intersection $-1$ and one can blow them down to obtain a projective plane with six distinct points $p_1,\ldots,p_4,d_5,d_6$ (abusing notation, we denote with the same letter points and the exceptional divisors above them). The five exceptional divisors of the blowing up of the quadric are, in the blowing up of the plane at the six points $p_1,p_2,p_3,p_4,d_5,d_6$, the exceptional divisors above $p_1,p_2,p_3,p_4$ and the strict transform of the line $(d_5d_6)$. The $S_5$ so obtained has to transversals, which are the strict transforms of two conics in the plane through $p_1,p_2,p_3,p_4,d_5$ and $p_1,p_2,p_3,p_4,d_6$, respectively.
\quad$\square$\par\medskip

\begin{figure}[h]\label{quadplan}
\caption{Blowing up the quadric vs blowing up the plane}
\begin{center}
\includegraphics[scale=.7]{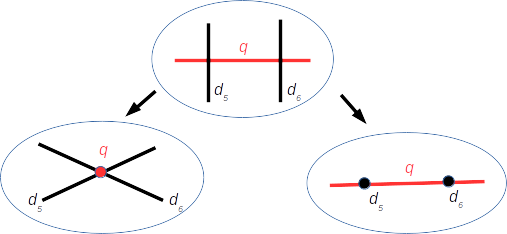}
\end{center}
\end{figure}

\subsection{Case study, the end}\label{Case3}

The part at infinity of the cubic surface, which is $\Sing_R^{\infty}=\Sing_R^h\setminus\Sing_R$, is the nonsingular cubic curve in the projective plane with homogeneous coordinates $x,y,z$ given by the equation
\begin{equation}\label{Singinfcase}
80x^3-107yx^2-47zx^2-9y^2x+95zxy-96z^2x-68y^3+98zy^2+78z^2y-24z^3
\end{equation}
The formulas (\ref{reciscrew}) for the reciprocal twist extend to the following homogeneous quadratic formulas:
\begin{equation}\label{RecTwinfCase}
\begin{aligned}
\Omega_{\infty}&= (-53\,x+92\,y-108\,z)\,\begin{pmatrix} x\\ y\\ z\end{pmatrix}\\
V_\infty&=\begin{pmatrix}
-80\,x^2+55\,x\,y-109\,x\,z+140\,y^2-14\,y\,z-84\,z^2\\
52\,x^2-131\,x\,y+164\,x\,z+68\,y^2-62\,y\,z+24\,z^2\\
156\,x^2-245\,x\,y+180\,x\,z-36\,y^2-102\,y\,z+24\,z^2
\end{pmatrix}
\end{aligned}
\end{equation}
These formulas have two points of indetermination on $\Sing_R^{\infty}$ which are, of course, the points at infinity of the two lines of indetermination for the formulas (\ref{reciscrew}) on $\Sing_R$ corresponding to $F_2$.

\begin{figure}[h]\label{IndetinfPts}
\caption{Points of indetermination at infinity of the quadratic formulas for the reciprocal twist}
\begin{center}
\includegraphics[scale=.5]{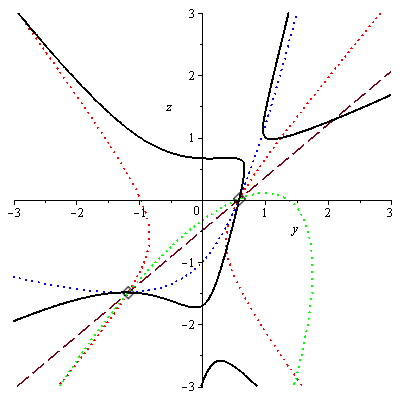}
\end{center}
\end{figure}

The figure \ref{IndetinfPts} shows the two points of indetermination in the affine part of the plane at infinity given by $x=1$. The solid black curve is the cubic $\Sing_R^{\infty}$, the dashed brown line is the line whose equation is the linear form factor of $\Omega_{\infty}$, and the dotted conics in red, green and blue are those whose equations are the components of $V_{\infty}$. All these curves concur to the indetermination points marked with diamonds.

We now check that the reciprocal twist mapping is also defined at the indetermination points by other formulas: precisely, we use the degree $3$ homogeneous part of the cofactors of the coefficients in the last line of the Jacobian matrix $\Jac$:

\begin{equation}\label{RecTwinfCubic}
\begin{aligned}
\Omega_{\infty}&= (-x^2+7yx-12zx-6z^2+7y^2-7zy)\,\begin{pmatrix} x\\ y\\ z\end{pmatrix}\\
V_\infty&=\begin{pmatrix}
-4 y x^2-12 z x^2+17 y^2 x+z x y-18 z^2 x+4 y^3+7 z y^2-5 z^2 y-6 z^3\\
4 x^3-17 y x^2+14 z x^2-4 y^2 x+11 z x y+6 z^2 x+2 z y^2+4 z^2 y\\
12 x^3-15 y x^2+18 z x^2-18 y^2 x-z x y+6 z^2 x-2 y^3-4 z y^2
\end{pmatrix}
\end{aligned}
\end{equation}

The homogeneous ideal generated by the ideal of indetermination of the quadratic formulas (\ref{RecTwinfCase}) and the quadratic form $-x^2+7yx-12zx-6z^2+7y^2-7zy$ factor of $\Omega_{\infty}$ in (\ref{RecTwinfCubic}) contains the third power of the maximal ideal generated by $x,y,z$. This shows that the reciprocal twist mapping is well defined as a regular map on the whole of $\Sing_R^h$.

\section*{Conclusion}

We have proved the rationality of the locus of singular configurations of the general Gough-Stewart platform. We have moreover related the rationality with the reciprocal twist mapping, which has a kinematic relevance; the lines of reciprocal twists form a quadric surface, and this fact may be useful for further studies on the singularities of a Gough-Stewart platform. We have also related a group of five lines on the cubic surface of singularities with special singular configurations for which the reciprocal twist is the twist of a pure rotation. This interplay between the classical algebraic geometry of cubic surfaces and kinematic properties of a parallel robot is rather fascinating.

We conclude with two questions encountered in this paper for which we have no satisfactory answer.

The first question concerns the characterization of singularities at infinity in section \ref{Secsinginf}. We have seen that it can be expressed in terms of a relation between two planar hexagons (Equation \ref{hexagons}). We have not been able to uncover the geometric significance of this relation. In the case of the planar 3-R\underline{P}R, the singularities at infinity may be characterized by the fact that two triples of aligned points can be transformed one into another by an affine mapping \cite{Co}. 

The second question is related to the fact that the reciprocal twist mapping extends to a regular mapping defined on the whole projective closure of the surface of singularities. We have followed a rather cumbersome way, with the help of a computation in a specific example, to show that this is indeed generically the case. It would be much nicer if the following assertion were true.
\begin{quote}
\textit{Let $A(t)\,X=0$ be a homogeneous linear system of $n$ equations in $n$ unknowns with coeffcients $a_{i,j}(t)$ polynomials of degree $1$ in parameters $t=(t_1,\ldots,t_p)$. Assume that the projective closure $S^h$ (in the $p$-dimensional projective space) of the set $S$ of parameters $t$ such that $\det(A(t))=0$ is a smooth projective hypersurface. Then the mapping which associates to $t\in S$ the line of solutions of $A(t)\,X=0$ extends to a regular mapping from $S^h$ to the $n-1$-dimensional projective space. }
\end{quote}
This assertion is not obvious when the degree of $\det(A(t))$ is strictly smaller than $n$, which happens in our case since the degree of the equation of the surface of singularities w.r.t. the position variables is $3$ instead of $6$. We have no idea whether this assertion holds true.

\bibliographystyle{plain}
\bibliography{referencesGS}
\end{document}